\documentclass[12pt]{article}

\usepackage{graphicx}
\usepackage{amsmath,amsthm,amssymb,color}
\usepackage[noadjust,sort]{cite}

\date{}

\renewcommand{\baselinestretch}{1.1}\normalsize

\numberwithin{equation}{section}

\def\({\bigl(}
\def\){\bigr)}

\def\abs#1{\lvert#1\rvert} 
 
\def\dfrac#1#2{\lower0.15ex\hbox{\large$\textstyle\frac{#1}{#2}$}}
\def\({\bigl(}
\def\){\bigr)}
\def\st{\mathrel{|}}
\def\Aut{\operatorname{Aut}}

\def\R{\mathcal{R}}
\def\I{\mathcal{I}}

\def\nicebreak{\vskip 0pt plus 50pt\penalty-300\vskip 0pt plus -50pt }

\begin{document}

\title{A class of Ramsey-extremal hypergraphs}

\author{Brendan~D.~McKay\vrule width0pt height2ex\thanks
 {Research supported by the Australian Research Council.}\\
\small Research School of Computer Science\\[-0.9ex]
\small Australian National University\\[-0.9ex]
\small Canberra ACT 2601, Australia\\[0.3ex]
\small\texttt{brendan.mckay@anu.edu.au}
}

\maketitle

\begin{abstract}
 In 1991, McKay and Radziszowski proved that, however each
 3-subset of a 13-set is assigned one of two colours, there
 is some 4-subset whose four 3-subsets have the same colour.
 More than 25 years later, this remains the only non-trivial
 classical Ramsey number known for hypergraphs.
 In this article, we find all the extremal colourings of the
 3-subsets of a 12-set and list some of their properties.
 Using the catalogue, we provide an answer to a question of
 Dudek, Fleur, Mubayi and R\"odl about the size-Ramsey numbers of
 hypergraphs.
\end{abstract}

\section{Introduction}\label{intro}
A colouring of all the $s$-subsets of an $n$-set with two 
colours is called \textit{$R(j,k;s)$-good} if there is no $j$-subset
(of the $n$-set) containing only $s$-subsets of the first
colour, and no $k$-subset containing only $s$-subsets of the
second colour. (Note that it is the $s$-subsets receiving colours,
not the elements of the $n$-set.)
The \textit{Ramsey number} $R(j,k;s)$ is
defined to be the least $n$ for which there is no
$R(j,k;s)$-good colouring.

Although there are several known values of $R(j,k;2)$~\cite{SRN}, 
which is usually written as just $R(j,k)$, the only known
non-trivial value of $R(j,k;s)$ for $s\ge 3$ is $R(4,4;3)=13$.
As a lower bound, a suitable colouring of the 3-subsets of
a 12-set was presented by Isbell in 1969~\cite{Isbell}, and this was
proved best possible by the present author and Radziszowski
in 1991~\cite{MR}.
During that project we found more than 200,000 $R(4,4;3)$-good
colourings for 12 points, but did not have the resources to
compute them all. With the aid of an improved algorithm
and the much greater computing resources available today, we
can now show that the number of $R(4,4;3)$-good colourings for
12 points is precisely 434,714.
We hope that this compilation of data will assist further
investigations.

\nicebreak
\section{Method}\label{method}
We prefer to use slightly different terminology for this
description.  Suppose we have an $R(4,4;3)$-good colouring
of the 3-subsets of an $n$-set $V$.
We will call the 4-subsets of $V$ \textit{quadruples}.

If we choose just the 3-subsets of $V$ having the first colour,
we obtain a 3-uniform hypergraph on $V$ with the property
that every quadruple contains 1, 2 or 3 edges (the other
possibilities 0 and 4 being forbidden). We will call this
a $R(4,4;3)$-good hypergraph.
Note that we could have chosen the other colour instead and
would have obtained the complementary hypergraph.
We can obvious recover the colouring from the hypergraph, so
we lose nothing by continuing with hypergraph terminology.

Denote by $\R(n)$ the set of $R(4,4;3)$-good hypergraphs
with $n$ points.  If we wish to emphasize the point set $V$,
we may write $\R(V)$ instead. More generally, $\R(n,e)$ is the
set of $R(4,4;3)$-good hypergraphs with $n$ points and $e$ edges,
and notations like $\R(V,{\le}110)$ have their obvious meanings.

Our aim is to find $\R(12)$.  By the remark just made, it
will suffice to find $\R(12,{\le}110)$, since
$110=\frac12\binom{12}{3}$ and the rest are complements.
Given $G\in\R(V)$ and $v\in V$, define $G_v$ to be the hypergraph
with point set $V{-}v$ and all the edges of $G$ that lie in~$V{-}v$.
Clearly $G_v\in\R(V{-}v)$.
Since the points of $G\in\R(n,e)$ lie on average in $3e/n$
edges, we find that for $G\in R(12,{\le}110)$ there is some
$v$ such that $G_v\in R(11,{\le}82)$.
Continuing such logic we find a construction path
\begin{equation}\label{path}
  \R(9,{\le} 41) \to \R(10,{\le} 59) \to \R(11,{\le} 82) 
  \to \R(12,{\le} 110).
\end{equation}
Each step in~\eqref{path} involves adding one point and some
edges that include the new point.  Moreover, we can assume that
the new point is in at least as many edges as any of the old points
(after the new edges are added).

The programs developed for~\cite{MR} are fast enough to find
$\R(9,{\le} 41)$ in a few hours.  There are exactly
3,030,480,232 such hypergraphs and these form our starting point.
It would be convenient to perform each of the three steps of~\eqref{path}
separately, but it would be quite expensive.  The number of hypergraphs in
$\R(10)$ and $\R(11)$ is greater than $10^{11}$ and
even the task of extending one hypergraph by one point requires
solution of  a large set of integer inequalities. We need a better way.

If $S$ is a set and $B\subseteq T\subseteq S$, then the
\textit{interval} $[B,T]$ is
$\{ X\subseteq S \st B\subseteq X\subseteq T\}$.  The use of
intervals for solving sets of inequalities efficiently
was introduced in~\cite{R45}.

Define $V_9=\{0,1,\ldots,8\}$ and $V_{10}=V_9\cup\{a\}$.
Consider extending $G_9\in \R(V_9)$ to all possible
$G_{10}\in \R(V_{10})$ by adding the point~$a$ and some edges
that include~$a$.  The possible edges all have the form
$\{i,j,a\}$ where $i,j\in V_9$; number these $e_0,e_1,\ldots,e_{35}$
in some order.  Each solution for $G_{10}$ corresponds to a subset
of $\{e_0,e_1,\ldots,e_{35}\}$.

Now consider the constraints required for $G_{10}$ to be $R(4,4;3)$-good.
The quadruples within $V_9$ are fine already, since we are not adding
any further edges inside $V_9$. So consider a quadruple 
$\{i,j,k,a\}$, where $i,j,k\in V_9$.  If $\{i,j,k\}$ is an edge of $G_9$,
we need that at least one of the edges $\{i,j,a\}, \{i,k,a\},\{j,k,a\}$ 
is not selected, while if $\{i,j,k\}$ is not an edge of $G_9$, at least
one of those three edges must be selected.

Now we can describe how intervals are used to process many cases
simultaneously. Consider one interval
$[B,T]\subseteq \{e_0,e_1,\ldots,e_{35}\}$
and one quadruple $\{i,j,k,a\}$. Define $X=\{e_r,e_s,e_t\}$,
where $e_r=\{i,j,a\}$, $e_s=\{i,k,a\}$ and $e_t=\{j,k,a\}$.
Now we apply the following \textit{collapsing rules}:

\def\strut{\vrule width0pt height2.4ex}
\begin{figure}[ht]
\centering
\begin{tabular}{c|c|c|l}
 $\{i,j,k\}\in G\,$? & $B\cap X$ & $T\cap X$ & replace $[B,T]$ by \\
 \hline \strut
 NO & $\ne\emptyset$ & any & $[B,T]$ \\
    & $\emptyset$    & $\emptyset$ & nothing \\
    & $\emptyset$    & $\{i\}$ & $[B{+}i,T]$ \\
    & $\emptyset$    & $\{i,j\}$ & $[B{+}i,T], [B{+}j,T{-}i]$ \\
    & $\emptyset$    & $\{i,j,k\}$ & $[B{+}i,T], [B{+}j,T{-}i], [B{+}k,T{-}i{-}j]$ \\
\hline\hline
 \strut$\{i,j,k\}\in G\,$? & $\bar T\cap X$  & $\bar B\cap X$ & replace $[B,T]$ by \\
\hline \strut
YES & $\ne\emptyset$ & any & $[B,T]$ \\
    & $\emptyset$    & $\emptyset$ & nothing \\
    & $\emptyset$    & $\{i\}$ & $[B,T{-}i]$ \\
    & $\emptyset$    & $\{i,j\}$ & $[B,T{-}i], [B{+}i,T{-}j]$ \\
    & $\emptyset$    & $\{i,j,k\}$ & $[B,T{-}i], [B{+}i,T{-}j], [B{+}i{+}j,T{-}k]$ \\
\end{tabular}
\caption{Collapsing rules for an interval $[B,T]$ based on quadruple $\{i,j,k,a\}$.
\label{fig}}
\end{figure}

By considering each case, we find that the effect of the collapsing
rules is to replace $[B,T]$ by a set of disjoint intervals whose
union is the set of all sets in $[B,T]$ that satisfy the quadruple
$\{i,j,k,a\}$.
For best practical performance, subsets of $\{e_0,e_1,\ldots,e_{35}\}$
can be represented by the bits in a single machine word, then
the collapsing rules can be implemented in a few machine
instructions each.

Starting with the interval $[\emptyset,\{e_0,e_1,\ldots,e_{35}\}]$
we apply the collapsing rules for each quadruple $\{i,j,k,a\}$.
The result is a set of disjoint intervals (typically a few hundred)
whose union gives exactly the set of all extensions of $G_9$ to $\R(10)$.
The efficiency depends a lot on the order in which quadruples are
processed; we found a good order by experiment.

Now consider further extension to $R(4,4;3)$-good hypergraphs
on $V_{11}=\{0,\ldots,8,a,b\}$.
The edges we need to add in total to $G_9$ either have the form
$\{i,j,a\}$ (already added in making $G_{10}$), 
$\{i,j,b\}$, or $\{i,a,b\}$, where in each case $i,j,k\in V_9$.
Here we can make an observation that is key to the whole computation:
\textit{The sets of edges $\{i,j,b\}$ which satisfy quadruples of the
form $\{i,j,k,b\}$ are the same as the sets of edges $\{i,j,a\}$
which satisfy quadruples of the form $\{i,j,k,b\}$, except that 
$a$ is replaced by~$b$.}

Given this observation, we make the possibilities for $G_{11}$
as follows, given $G_9$, a set $\I$ of intervals describing the
extensions of $G_9$ to $\R(10)$, and a particular extension~$G_{10}$.
The possible new edges are numbered $e_0,\ldots,e_{44}$, where
$e_0,\ldots,e_{35}$ are edges of the form $\{i,j,b\}$ numbered in
the same order as we numbered the edges $\{i,j,a\}$ in the previous
step, and $e_{36},\ldots,e_{44}$ are the edges of the form
$\{i,a,b\}$ in any order. To find all solutions, instead of
starting with the single interval
$[\emptyset,\{e_0,e_1,\ldots,e_{44}\}]$ as in the previous step,
we start with the set of intervals
$[B,T\cup\{e_{36},\ldots,e_{44}\}]$ for $[B,T]\in\I$.  Then we avoid collapsing
rules which are unnecessary for the stated reasons.  This
results in a massive speedup.

To complete the process by extending from 11 to~12 points, we
use the same idea to begin with intervals obtained during the
extension to~11 points.  This phase is very fast as most intervals
are destroyed very quickly and only a comparatively small number
of solutions are found.

\medskip

It would be possible to apply the general method of~\cite{orderly}
to perform exhaustive isomorph reduction at each step in the
computation, but the large number of intermediate hypergraphs
makes that unwise.  Instead, we applied a weaker filter.
For a hypergraph with points $V$ and point $v\in V$, define
$d_v$ to be the number of edges that include~$v$.
Also define $f_v = \sum_e d_vd_wd_x$, where the sum is over
all edges $e=\{v,w,x\}$ that include~$v$.
Suppose we make $G\in G(V)$ by extending a smaller hypergraph,
and that $v\in V$ is the last point added.
The construction path~\eqref{path} assumed that $d_v\ge d_w$
for all $w\in V$, so that is the first filter applied.
If that doesn't eliminate $G$, we also require that $f_v$ be
maximum out of all $w\in V$ with maximum~$d_w$.  These rules
eliminate most isomorphs and are fast to apply.  When we finally
have a collection of $R(4,4;3)$-good hypergraphs on 12 points,
we perform complete isomorphism reduction using \texttt{nauty}~\cite{nauty}.

\begin{table}[p]
\centering
\begin{tabular}{c|c|c||c|c|c}
 ~~$n$~~ & ~~$e$~~ & count & ~~$n$~~ & ~~$e$~~ & count \\
 \hline
 3 & 0 & 1    &                                 9 & 33 & 2 \\
 \multicolumn{2}{c|}{total} & 2   &       & 34 & 204 \\
\cline{1-3}
 4 & 1 & 1    &                                    & 35 & 22616 \\
    & 2 & 1    &                                    & 36 & 774043 \\
 \multicolumn{2}{c|}{total} & 3   &       & 37 & 10877731 \\
 \cline{1-3}
 5 & 3 & 1    &                                    & 38 & 79336073 \\
    & 4 & 3    &                                    & 39 &  341024774 \\
    & 5 & 4    &                                    & 40 & 928650036 \\
 \multicolumn{2}{c|}{total} & 12  &      & 41 & 1669794753 \\
 \cline{1-3}
 6 & 6 & 1    &                                   & 42 & 2025923846 \\
    & 7 & 5    &      \multicolumn{2}{c|}{total} & ~~8086884310~~ \\
  \cline{4-6}
    & 8 & 22  &                              10 & 50 & 13 \\
    & 9 & 50  &                                   & 51 & 1810 \\
    & 10 & 70 &                                  & 52 & 121356 \\
 \multicolumn{2}{c|}{total} & 226 & \multicolumn{2}{c|}{$\cdots$}\\
 \cline{1-3}
 7 & 12 & 1 &        \multicolumn{2}{c|}{total} & $\approx 6.2{\times}10^{11}$ \\
 \cline{4-6}
    & 13 & 26 &                             11 & 73 & 36 \\
    & 14 & 338 &                                & 74 & 4725 \\
    & 15 & 1793 &                              & 75 & 246299 \\
    & 16 & 5055 &                   \multicolumn{2}{c|}{$\cdots$} \\
    & 17 & 8317 &   \multicolumn{2}{c|}{total} & $\approx 2.1{\times}10^{11}$ \\
    \cline{4-6}
  \multicolumn{2}{c|}{total} & 31060 &   12 & 104 & 4 \\
    \cline{1-3}
 8 & 21 & 1 &                                             & 105 & 123 \\
    & 22 & 278 &                                         & 106 & 1465 \\
    & 23 & 9763 &                                       & 107 & 10235 \\
    & 24 & 107241 &                                   & 108 & 41939 \\
    & 25 & 573596 &                                   & 109 & 98235 \\
    & 26 & 1764747 &                                 & 110 & 130712 \\
    & 27 & 3380337 &  \multicolumn{2}{c|}{total} & 434714 \\
    & 28 & 4182459 \\
  \multicolumn{2}{c|}{total} & ~~15854385~~
\end{tabular}
\caption{The numbers of $R(4,4;3)$-good hypergraphs with $n$ points and
  $e$ edges.  The totals include complements. \label{counts}}
\end{table}

\nicebreak
\section{Results}\label{results}

There are about $8.4\times 10^{11}$ $R(4,4;3)$-good hypergraphs
altogether, including 434,714 with 12 points.
Table~\ref{counts} details the numbers of $R(4,4;3)$-good hypergraphs
for each number of points and edges.
For 10 and 11 points we only did incomplete isomorph reduction, as
explained above; hence the totals for those sizes are estimates.

The \textit{automorphism group} $\Aut(G)$ of a hypergraph $G\in\R(V)$
is the set of permutations of $V$ which preserve the edge set.  As detailed
in Table~\ref{groups}, most hypergraphs in $\R(12)$ have a trivial group
and none have a transitive group.  The unique hypergraph with
$\abs{\Aut(G)}=60$, which has two orbits of size~6, is presented in
Figure~\ref{example} using letters for elements of~$V$.
This hypergraph is one of the 1306 in $\R(12)$ that are
self-complementary and is isomorphic to the one found by 
Isbell~\cite{Isbell}.

\begin{table}[th]
\centering
\begin{tabular}{c|c|c}
  $\abs{\Aut(G)}$ & orbits & count \\
  \hline
  1 & 12 & 432300 \\
  2 & 6 & 18 \\
     & 7 & 112 \\
     & 8 & 1669 \\
  3 & 4 & 529 \\
  4 & 6 & 32 \\
  6 & 2 & 20 \\
     & 4 & 17 \\
10 & 4 & 1 \\
12 & 2 & 15 \\
60 & 2 & 1 \\
\end{tabular}
\caption{Counts of $\R(12)$ by automophism group.\label{groups}}
\end{table}

None of the hypergraphs in $\R(12)$ extend to a hypergraph in
$\R(13)$, consistently with the finding of~\cite{MR} that
$\R(13)=\emptyset$.  This raises the question of how close we can
get to a hypergraph in $\R(13)$; specifically, how many edges
of the complete hypergraph $K_{13}^{(3)}$ can we colour without
obtaining a monochromatic induced $K_4^{(3)}$?
The generation method described in the previous section
can be easily adapted to ignore particular quadruples.  If we ignore
the constraints normally attributed to the quadruples $\{i,j,k,a\}$ which
contain a specified $\{i,j,a\}$, then we are colouring the edges of
the complete hypergraph except for one uncoloured edge.
Using this method we found that $K_{13}^{(3)}$ minus one edge
cannot be coloured with two colours
without creating a monochromatic $K_4^{(3)}$.

On the other hand, if we omit two edges of $K_{13}^{(3)}$, a colouring
without a monochromatic $K_4^{(3)}$ may be possible. In
Figure~\ref{minus2} we give examples where the two omitted
edges overlap in one or two points.  We did not find any examples
with the omitted two edges being disjoint, but our search in that case
was not exhaustive.  We can report these partial results: there is no
good colouring of $K_{13}^{(3)}$ minus the edges $\{1,2,3\}$ and
$\{4,5,6\}$ such that a good colouring of $K_{12}^{(3)}$ can be
obtained either by
deleting vertex 1 and colouring edge $\{4,5,6\}$, or by deleting
vertex 7 and colouring both edges $\{1,2,3\}$ and $\{4,5,6\}$.
We propose the remaining cases of two disjoint edges as a
challenge for the reader.

If $H$ is a 3-uniform hypergraph, the \textit{size-Ramsey number}
$\hat R^{(3)}(H)$ is the least number $m$ such that for some
3-uniform hypergraph $G$ with $m$ edges, every colouring
of the edges of $G$ with two colours includes a monochromatic copy
of~$H$.  If $H=K_4^{(3)}$, then the value $R(4,4;3)=13$ implies that 
$\hat R^{(3)}(H)\le\binom{13}{3}=286$ since we can take
$G=K_{13}^{(3)}$.
Dudek, La~Fleur, Mubayi and R\"odl~\cite[Question\,2.2]{Dudek}
ask whether this bound is sharp.
Since $K_{13}^{(3)}$ minus one edge cannot be coloured without creating
a monochromatic $K_4^{(3)}$ we have $\hat R^{(3)}(H)\le 285$, which
answers Dudek et~al.'s question in the negative.

The extremal $R(4,4;3)$-good hypergraphs are available
online~\cite{online}.
Finally, we thank Staszek Radziszowski for many useful
comments.

\begin{figure}[ht]
\centering
\medskip
\parbox{0.74\hsize}{\small
Let
$\varGamma=\langle(\texttt{cd})(\texttt{ef})(\texttt{CD})(\texttt{EF}),
 \allowbreak
        (\texttt{bc})(\texttt{de})(\texttt{BC})(\texttt{DE}),
 \allowbreak
         (\texttt{ab})(\texttt{ef})(\texttt{AB})(\texttt{EF})\rangle$
be a permutation acting on the points \texttt{abcdefABCDEF}.
It is isomorphic to the alternating group $A_5$ and
acts 2-transitively on
each of its orbits $\{\mathtt a,\ldots,\mathtt f\}$
and $\{\mathtt A,\ldots,\mathtt F\}$.
Now construct a hypergraph by applying $\varGamma$ to each
of the starting edges $\{\texttt{abe}, \texttt{ABE}, \texttt{abC},
 \texttt{aAB}, \texttt{cAB}\}$.
These provide 10, 10, 30, 30 and 30 edges, respectively.
The hypergraph induced by each orbit is the same
2-(6,3,2) design. The relabelling
\texttt{(aD)(bC)(cB)(dA)(eF)(fE)} takes the hypergraph
onto its complement.}
\caption{The unique hypergraph in $\R(12)$ with
automorphism group of order 60.\label{example}}
\end{figure}

\begin{figure}[ht]
\centering
\medskip
\parbox{0.43\hsize}{\small\baselineskip=12pt\tt
 acd bcd abe ace bce cde adf cdf \\
 def adg aeg beg ceg deg afg bfg \\
 efg ach bch adh bdh aeh beh deh \\
 afh efh bgh dgh fgh bdi cdi bei \\
 cei afi bfi dfi efi bgi cgi dgi \\
 ahi chi dhi ehi abj acj bcj cdj \\
 aej dej bfj cfj agj bgj cgj dgj \\
 bhj ehj fhj aij gij abk bck bdk \\
 cek afk bfk cfk efk cgk fgk dhk \\
 ghk aik bik eik hik bjk djk gjk \\
 hjk ijk abl acl bcl adl bel afl \\
 bfl cfl bgl dgl chl fhl ghl eil \\
 fil gil hil ajl djl ejl fjl ijl \\
 akl ckl dkl ekl hkl abm adm bdm \\
 aem dem bfm cfm dfm cgm egm ahm \\
 bhm chm ghm aim cim fim ejm fjm \\
 hjm ijm akm dkm ekm fkm gkm jkm \\
 blm clm dlm elm glm ilm \\
{\rm Omitted edges:} abc ade}
\parbox{0.43\hsize}{\small\baselineskip=12pt\tt
 bcd cde acf bcf aef def adg bdg \\
 cdg aeg beg deg bfg efg abh ach \\
 bch adh bdh beh cfh egh fgh aci \\
 aei bei cei dei afi dfi agi cgi \\
 fgi bhi dhi fhi adj bdj aej cej \\
 bfj cfj dfj agj bgj ahj bhj chj \\
 ehj ghj bij dij eij abk ack bck \\
 cdk aek bek cek cfk dfk efk bgk \\
 cgk fgk ahk dhk fhk aik bik gik \\
 hik ejk fjk gjk hjk ijk abl acl \\
 bcl adl bdl afl bfl dfl efl cgl \\
 dgl fgl chl dhl ehl ghl bil cil \\
 eil fil hil ajl cjl ejl fjl gjl \\
 bkl dkl gkl abm adm bdm cdm bem \\
 cem afm bfm cfm dfm efm agm bgm \\
 ahm ehm fhm ghm cim fim gim cjm \\
 djm gjm ijm dkm ekm hkm jkm blm \\
 elm ilm jlm klm \\
{\rm Omitted edges:} abc abd}
\caption{Two $R(4,4;3)$-good colourings of the complete hypergraph 
$K_{13}^{(3)}$ minus two edges. Edges not mentioned have
the second colour.\label{minus2}}
\end{figure}

\nicebreak
\renewcommand{\baselinestretch}{1}

\end{document}